\documentclass[a4paper,12pt,reqno]{amsart}
\usepackage{rotating,pdflscape,color,amssymb,subfigure,psfrag,amsmath,eufrak,bbm,epsfig}
\usepackage{a4wide,amscd}
\usepackage{tikz} % doc : chercher PGF
 \usepackage[usenames,dvipsnames]{pstricks}
 \usepackage{pst-grad} % For gradients
 \usepackage{pst-plot} % For axes

\definecolor{vdarkred}{rgb}{0.6,0,0.2}
\definecolor{vdarkblue}{rgb}{0,0.2,0.6}
\usepackage[pdftex, colorlinks, linkcolor=vdarkblue,citecolor=vdarkred]{hyperref}
\usepackage[small]{caption}

\newcommand{\Cc}{\mc{C}}
\newcommand{\Ga}{\Gamma} 
\newcommand{\ga}{\gamma}

\newcommand{\ld}{\ldots}
\newcommand{\beg}{\begin}
\newcommand{\en}{\end}

\newcommand{\cdd}{\cdots\cdots}
 
\newcommand{\trm}{\textrm}
 
\newcommand{\bgt}{\begin{itemize}}
\newcommand{\ent}{\end{itemize}}
\newcommand{\ite}{\item}

\newcommand{\op}{\operatorname}
\newcommand{\eqre}{\eqref}
\newcommand{\re}{\ref}
\newcommand{\la}{\label}

\newcommand{\si}{\sigma}

\newcommand{\Var}{\operatorname{Var}}

\newcommand{\lan}{\langle}
\newcommand{\ran}{\rangle}

\newcommand{\ds}{\displaystyle}
 
\newcommand{\p}{\mathbb{P}}

\newcommand{\Tr}{\operatorname{Tr}}

\newcommand{\ninf}{\underset{N\to\infty}{\longrightarrow}}

\newcommand{\Ninf}{\underset{N\to\infty}{\longrightarrow}}

\newcommand{\E}{\op{\mathbb{E}}}

\newcommand{\R}{\mathbb{R}}
\newcommand{\C}{\mathbb{C}}

\newcommand{\z}{\mathbb{Z}}

\newcommand{\ud}{\mathrm{d}}

\newcommand{\pro}{probability }

\newcommand{\f}{\frac}
\newcommand{\ff}{\frac{1}}
\newcommand{\lf}{\left}
\newcommand{\ri}{\right}

\newcommand{\st}{such that }

\newcommand{\lam}{\lambda}

\newcommand{\ti}{\times}

\newcommand{\vfi}{\varphi}
\newcommand{\ste}{\, ;\, }
\newcommand{\mc}{\mathcal }

\newcommand{\eps}{\varepsilon}

\newcommand{\bck}{\backslash}
\newcommand{\al}{\alpha}
\newcommand{\tta}{\theta}

\newcommand{\ovl}{\overline}

\newcommand{\bbm}{\begin{bmatrix}}
\newcommand{\ebm}{\end{bmatrix}}
\newcommand{\bes}{\begin{equation*}}
\newcommand{\ees}{\end{equation*}}
\newcommand{\be}{\begin{equation}}
\newcommand{\ee}{\end{equation}}
\newcommand{\beqy}{\begin{eqnarray}}
\newcommand{\eeqy}{\end{eqnarray}}
\newcommand{\beq}{\begin{eqnarray*}}
\newcommand{\eeq}{\end{eqnarray*}}
\newcommand{\one}{\mathbbm{1}}
\newcommand{\lto}{\longrightarrow}

\newcommand{\ie}{\emph{i.e. }}

\newcommand{\bpm}{\begin{pmatrix}}
\newcommand{\epm}{\end{pmatrix}}

\newcommand{\cd}{\cdots}

\newcommand{\X}{\mc{X}}

\newcommand{\hmu}{\hat{\mu}}

\newcommand{\bpr}{\beg{proof}}
\newcommand{\epr}{\en{proof}}

\newcommand{\Si}{\Sigma}
\newcommand{\bet}{\beta}
\newcommand{\del}{\delta}
\newcommand{\Del}{\Delta}

\newcommand{\otd}{^{\otimes 2}}

\newcommand{\nober}{\nonumber}

\newcommand{\pa}{\partial}

\newcommand{\tR}{\widetilde{R}}

\newcommand{\Tta}{\Theta}

  \newcommand{\mre}{\mathrm{e}}
   \newcommand{\ii}{\mathrm{i}}

\theoremstyle{definition}

\long\def\symbolfootnote[#1]#2{\begingroup
\def\thefootnote{\fnsymbol{footnote}}\footnote[#1]{#2}\endgroup}

\title[Poisson statistics for matrix ensembles]{Poisson statistics for matrix ensembles at large temperature}


\author{Florent Benaych-Georges} \address{MAP 5, UMR CNRS 8145 - Universit\'e Paris Descartes, 45 rue des Saints-P\`eres 75270 Paris cedex~6, France.} \email{florent.benaych-georges@parisdescartes.fr}

\author{Sandrine P\'ech\'e} \address{LPMA, Universit\'e Paris Diderot, 5 rue Thomas Mann
75013 Paris.} 
\email{sandrine.peche@math.univ-paris-diderot.fr}

\date{\today}

\begin{document}
\maketitle
 
 \beg{abstract}In this article, we consider $\beta$-ensembles, \ie collections of particles with random positions on the real line having joint distribution  $$\frac{1}{Z_N(\beta)}|\Delta(\lambda)|^\beta e^{- \frac{N\beta}{4}\sum_{i=1}^N\lambda_i^2}\ud \lambda,$$ in the regime where $\beta\to 0$ as $N\to\infty$. We briefly describe the global regime and then consider the local regime. In the case where $N\beta$ stays bounded, we prove that the local eigenvalue statistics, in the vicinity of any real number, are asymptotically to those of a Poisson point process. In the case where $N\beta\to\infty$, we prove a partial result in this direction.  
 \end{abstract}
 %arxiv abstract:
 %In this article, we consider $\beta$-ensembles, i.e. collections of particles with random positions on the real line having joint distribution  $$\frac{1}{Z_N(\beta)}|\Delta(\lambda)|^\beta e^{- \frac{N\beta}{4}\sum_{i=1}^N\lambda_i^2}d \lambda,$$ in the regime where $\beta\to 0$ as $N\to\infty$. We briefly describe the global regime and then consider the local regime. In the case where $N\beta$ stays bounded, we prove that the local eigenvalue statistics, in the vicinity of any real number, are asymptotically to those of a Poisson point process. In the case where $N\beta\to\infty$, we prove a partial result in this direction.
 %\tableofcontents
 
 \section{Introduction}
 General $\beta$-ensembles are collections of particles with random positions on the real line with joint distribution  
\be\la{1521216h392512}\ff{Z_N(\bet)}|\Del(\lam)|^\bet e^{- \sum_{i=1}^NV({\lam_i})}\ud \lam,\ee where $\ud\lam$ denotes the Lebesgue measure on $\R^N$, $\Del(\lam):=\prod_{1\le i<j\le N}(\lam_j-\lam_i)$, $V$ is a potential with enough growth at infinity (like the Gaussian potential $V_G(x)=\f{x^2}{2}$) and $Z_N(\bet)$ is a normalizing constant. 
Their study is  initially motivated by some considerations from physics:   
the probability distribution can be viewed as the equilibrium measure of a one dimensional Coulomb gas, but they  actually appear  to  be  connected  to  a  broad  spectrum  of  mathematics  and physics,   such as random matrices, number theory, lattice gas theory, quantum mechanics and Selberg-type integrals. 
In the case where  $\bet=1$, $2$ or $4$, the \pro measure of \eqre{1521216h392512} is the joint distribution of the eigenvalues of a random $N\ti N$ matrix $M$ with density proportional to $e^{- \Tr V(M)}$ on the space of respectively real symmetric, complex Hermitian or quaternionic Hermitian matrices (see e.g. \cite{alice-greg-ofer}). 
Besides, it was proved by Dumitriu and   Edelman in \cite{DumEdel} that when $V$ is  the Gaussian potential, for any $\bet>0$,  the \pro measure of \eqre{1521216h392512} is the joint distribution of the eigenvalues of the random $N\ti N$ tridiagonal matrix   
 \be\la{matH6513intro}H=  \bpm g_1 & X_{2}&&&\\
X_{2}& g_2& X_{3 }&&\\
&X_{3 }&\ddots&\ddots\\
&&\ddots&\ddots&X_{N}\\
&&&X_{N }&g_N
\epm,\ee
where the $g_i$'s are some $N(0,1)$ variables and for all $i$, $X_{i }=\sqrt{Y_{i}}$, with $Y_{i }$ distributed thanks to the $\Ga((i-1)\bet/2)$ law, everything being independent.  
For general potential $V$, there is   no random matrix representation.

In the classical cases ($\beta=1,2,4$), the combinatorial structure and repulsive interaction has been well--understood for a long time via the theory of determinantal or Pfaffian processes (see e.g. \cite{alice-greg-ofer} for references). 
The understanding of these asymptotic spectral statistics  to the full class of parameters $\beta>0$ has recently mobilized a lot of research. For general $\beta$, despite the lack of structure, some enormous progress has been accomplished recently.
For fixed $\beta$, a few results are now known. First, it is known  from \cite{BAG3} that the empirical eigenvalue distribution of the rescaled matrix 
$\frac{1}{\sqrt N}H $ converges weakly as $N \to \infty$ to a \pro measure which is the semi-circle distribution in the case of  Gaussian potential. 
The local eigenvalue statistics in the large $N$-limit are also quite well understood. In the Gaussian setting,
at the edge of the spectrum, Ram\'{\i}rez,  Rider and Vir\'ag have shown in \cite{RamirezRiderVirag} that the 
eigenvalues of $N^{1/6} (H -2\sqrt N I)$ converge in distribution  to those of the so-called stochastic Airy operator. 
In the bulk of the spectrum, the limiting spectral statistics are asymptotically defined in terms of the Sine-$\beta$ process, which is again defined as the solution of a stochastic equation by Valk\'o and Vir\'ag in  \cite{ValkoVirag}. In particular the authors show that the Sine-$\beta$, which is translation invariant, has a geometric description in terms of the Brownian carousel, a deterministic function of the Brownian motion in the hyperbolic plane. Some advances on $\bet$-ensembles have also been made by Sosoe and  Wang \cite{sosoe,sosoe2} and Bao and Su in \cite{bao}.\\
The question of universality for these statistics has now become an important matter of interest: some enormous progress has recently been accomplished by Bourgade,  Erd\"os and Yau  in \cite{BEY1,BEY2,BEY3,BEY4}. Therein the authors consider general $\beta$-ensembles (when the   potential   $V$   is $\mathcal{C}^4$ and regular,  or, in  the first papers,  convex and analytic). Assuming that the limiting spectral distribution (which   depends on $V$) is supported on a single interval, they prove that the limiting eigenvalue statistics at the edge of the spectrum are given by the $\beta$-Tracy-Widom distribution.  The universality in the bulk of the spectrum is also proved.\\
Another point of view to tackle $\beta$-ensembles and in particular the quantitative aspect of the repulsion between eigenvalues has been developed in particular by Allez, Bouchaud  and Guionnet in \cite{AllezGuionnet, AllezGuionnetBouchaud}. They show in particular that when $\beta \leq 2$, $\beta$-ensembles can be seen as    an $N$-dimensional process whose evolution is a mixing of that of $N$ independent real Brownian motions and of that of a  $\bet$-Dyson Brownian motion.

The scope of this article is to understand the  spectral behavior, at microscopic scale,  of $\beta$-ensembles in the case where $\beta\to 0$ and $N \to \infty$ (so that $\beta $  depends on the dimension $N$). At   macroscopic scale, such ensembles have been considered recently by \cite{TrinhTomoyuki} (see also the close model studied in \cite{Allez2}). Therein it is proved that when $\beta N \to c$ for some constant $c>0$, the scaled empirical eigenvalue distribution of $ \frac{1}{\sqrt \beta}H $ converges to the spectral measure of a deterministic Jacobi matrix, the density of which   is explicit.
When $\beta N \to \infty$, the limiting empirical eigenvalue distribution of $ $ converges  to the semi-circle distribution.
Local eigenvalue statistics have not been considered yet.

We here also consider the regime where  $\beta \to 0$ and $N \to \infty$, but study the local eigenvalue statistics. In \cite{KillipStoicu}, Killip and Stoiciu
 have considered the same question for circular $\beta$-ensembles. More precisely they study CMV matrices (which are discrete one-dimensional Dirac-type operators) with random decaying coefficients.  For rapidly decreasing coefficients, the eigenvalues have rigid spacing while in the case of slow decrease, the eigenvalues are distributed according to a Poisson process. More precisely,  they prove that local eigenvalue statistics of $\beta$-circular ensembles 
when $\beta\to 0$ are in the large $N$ limit those of a Poisson process.\\
For real-symmetric ensembles, the same question has recently been considered from a formal point of view.  Indeed,
in  \cite{AD1,AD2}, Allez and Dumaz    considered the $\beta \to 0$ limit  of the Sine-$\beta$ process and   of the $\beta$-Tracy-Widom  distribution. The $\beta\to 0$ limit of the Sine-$\bet$ is also considered by Lebl\'e and Sefaty in \cite{LebleSerfaty}. The approach used by \cite{LebleSerfaty} is based on approximation theory while \cite{AD1} use the diffusion representation of the Sine-$\beta$ process to consider the limit $\beta \to 0.$\\
One would expect again to prove that when $\beta \to 0$ simultaneously to $N \to \infty$, the eigenvalues in the vicinity of a point $u$ in the bulk of the spectrum exhibit Poisson statistics. In this text, we prove that this is true when $N\bet$ stays bounded as $N\to\infty$. In the case where  $\beta \to 0$ but $N\bet\to\infty$, we have a partial result which formally implies the Poisson statistics in the bulk, but does not allow to get a complete proof.\\ 
In Figure \re{Fig:space}, we compare this result with numerical simulations, giving a numerical evidence of the fact that the Poisson approximation works well (but gets less accurate as $\bet$ grows).
%we compare   histograms of spacings of the $\lam_k$'s from some numerical simulations with the density of the spacings in a Poisson point process (namely the exponential density) for several small values of   $\beta$. We see that the larger $\bet$ is, the less accurate the Poisson approximation seems to get. 
\begin{figure}[ht]
\centering
\includegraphics[scale=.4]{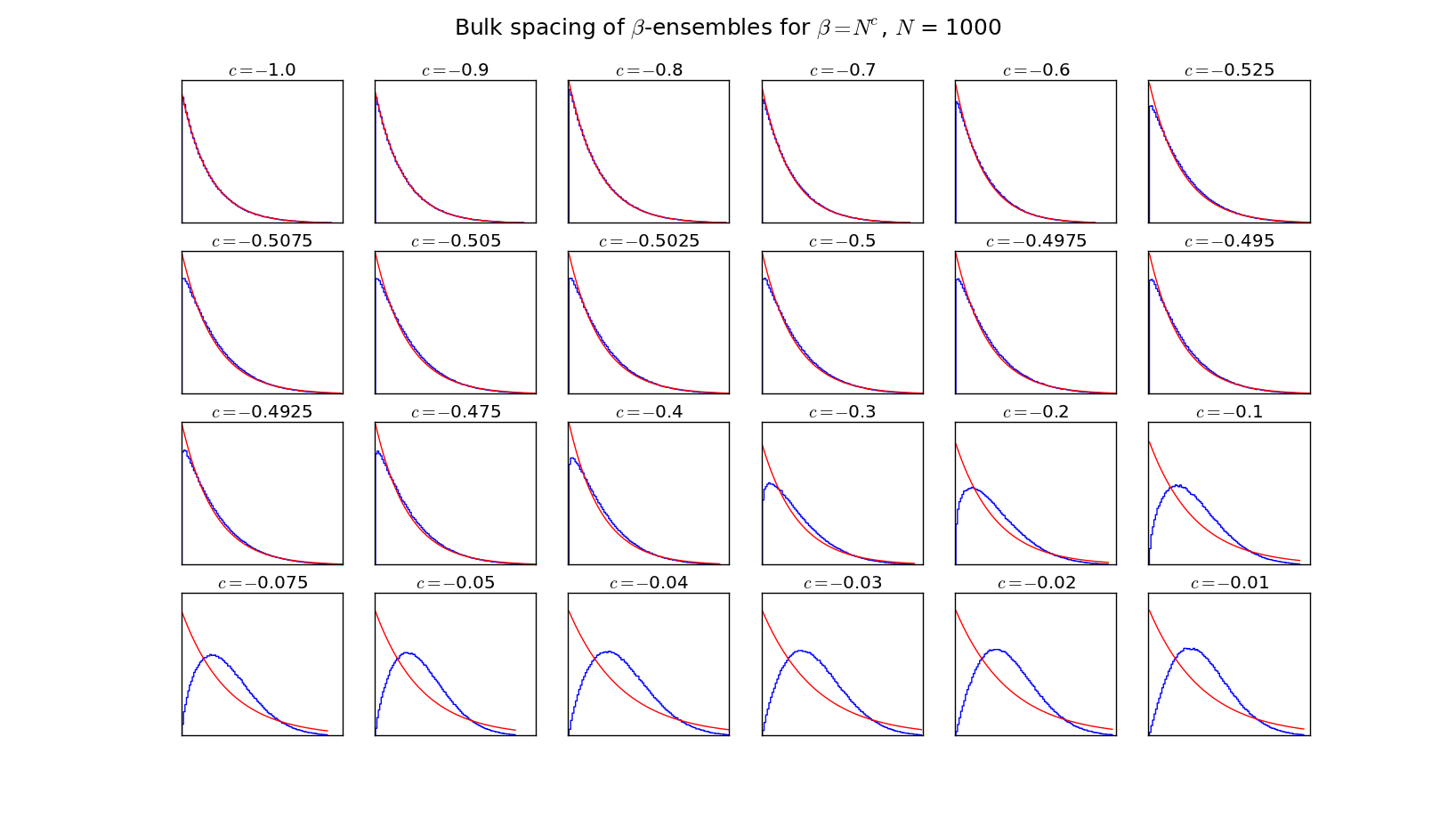}
\caption{{\bf Numerical comparison  of the spacings of the $\lam_k$'s    with the   spacings in a Poisson point process:} histogram, in blue,  of the debiazed spacings $(\lam_{k+1}-\lam_k)\sqrt{4-\lam_k^2}$ for $0.05N\le k\le .95N$, with $\beta=N^{c}$, for several values of $c\in (-1,0)$,  compared to the density, in red, of the  exponential law. The histogram of the spacings is drawn  thanks to $10^3$ independent   
realizations  of the distribution of \eqre{13313taokp} for $N=10^3$,     (we ordered the $\lam_k$'s in an increasing way so that $\lam_{k+1}-\lam_k$ is actually the spacing between two consecutive particles). We see that the Poisson approximation works well, but gets less accurate as $\bet$ grows  (note however that all simulations of the last row were made with $\bet$ almost equal to $1$).}\label{Fig:space}
\end{figure}

\paragraph{\bf Notation} For $u=u(N)$ and $v=v(N)$ some sequences, $$u\ll v\iff u/v\ninf 0\qquad ;\qquad u\sim v\iff u/v\ninf 1.$$

\paragraph{\bf Acknowledgment} \textit{The authors thank Alice Guionnet for her contribution to the proof and her useful suggestions for simplifying the arguments. We also thank Paul Bourgade for useful discussions.}

\newpage

\section{Statement of results}

 \subsection{Presentation of the model}

%\subsubsection{Particle model}
 For any $\al,\bet>0$ and any $N\ge 1$, we define \be\la{13313}Z_N(\al, \bet):=\int_{\lam\in \R^N} \Del(\lam)^\bet\mre^{-\f{\al}{2}\sum_{i=1}^N \lam_i^2}\ud \lam,\ee with $ \Del(\lam):=\prod_{1\le i<j\le N}|\lam_j-\lam_i|$.
 %For any $\al>0$ and any $N\ge 1$, we define \be\la{13313}Z_N(\al):=\int_{\lam\in \R^N} \Del(\lam)^\al \mre^{-\f{N\al}{4}\sum_{i=1}^N \lam_i^2}\ud \lam,\ee with $ \Del(\lam):=\prod_{1\le i<j\le N}|\lam_j-\lam_i|$.

 Let us  now  %consider $\bet=\bet(N)=2\ga N^{-1}$, and
  consider an exchangeable family  $(\lam_1, \ld, \lam_N)$ of random variables with    joint law     \be\la{13313taokp}P^{(N)}_{\al,\bet}(\ud \lam_1, \ld,\ud\lam_N ):=\f{1}{Z_N(\al,\bet)}\Del(\lam)^\bet \mre^{-\f{\al}{2}\sum_{i=1}^N \lam_i^2}\ud \lam_1\cd\ud\lam_N\ee
with $Z_N(\al,\bet)$ the normalization constant defined at \eqre{13313}.

\subsection{Tridiagonal model and relation between $\al$ and $\bet$}
Let %$\bet=2N^{-1}$ and  
 \be\la{matH6513}H= \ff{\sqrt{\al}}\bpm g_1 & X_{2}&&&\\
X_{2}& g_2& X_{3 }&&\\
&X_{3 }&\ddots&\ddots\\
&&\ddots&\ddots&X_{N}\\
&&&X_{N }&g_N
\epm,\ee
where the $g_i$'s are some $N(0,1)$ variables and for all $i$, $X_{i }=\sqrt{Y_{i}}$, with $Y_{i }$ distributed thanks to the $\Ga((i-1)\bet/2)$ law. We known, by \cite{DumEdel} or Section 4.5 of \cite{alice-greg-ofer}, that $P^{(N)}_{\al,\bet}$ is the joint law of the eigenvalues of $H$.

Note that $\Tr H$ is centered and $\al\E\Tr H^2=N+\bet\f{N(N-1)}{2}$, so that 
%to have $\al\E\Tr H^2\sim \al N$, one needs to have $1+\bet(N-1)/2\sim \al$. Hence if $N\bet\gg {1}$, we shall choose $\al\sim N\bet/2$, whereas if $N\bet\lto  2\ga\ge 0$ fixed, we shall choose $\al\lto \ga +1$.  A 
for \be\la{eq:elpha}\al\sim 1+N\bet/2,\ee
  the empirical eigenvalue distribution of $H$ has asymptotic first moments $0$ and $1$.

   \subsection{Global and local regime  for bounded $N\bet$}
 
 The following proposition gives the limit of the empirical  distribution of the $\lam_i$'s. The \pro measure $\mu_\ga$ in question here has been studied in \cite{AllezGuionnetBouchaud,TrinhTomoyuki}.
\beg{propo}[Global regime for bounded $N\bet$]\la{Icnbty413}Suppose that $N\bet \lto 2\ga\ge 0$ as $N\to\infty$ and $\al\lto \ga+1$. \bgt\ite[a)] 
Then under the law $P^{(N)}_{\al,\bet}$, the random \pro measure $$\ff{N}\sum_{i=1}^N \del_{\lam_i}$$ converges in \pro to an even \pro measure $\mu_\ga$  on $\R$, depending only on $\ga$,  with moments $m_k$ defined at \eqre{ictwyuuhb}, satisfying   $m_2=1$. \\
\ite[b)] The measure $\mu_\ga$  is absolutely continuous with respect to the Lebesgue measure on $\R$, with a density $f_{\mu_\ga}$ that is bounded on any compact set, and satisfies, for all $x>0$,    $$\mu_\ga(\R\bck[-x,x]) \;\le \; C_\ga\f{\mre^{- \f{\ga+1}{2}x^2 }}{x} $$where $C_\ga$ is a constant depending only on    $\ga$.\\
\ite[c)]   $\mu_\ga$ depends continuously on $\ga\ge 0$, is equal to  $N(0,1)$ if and only if $\ga=0$,  and   tends to the semicircle law with support $[-2,2]$ as $\ga\to\infty$.\\
%\ite[d)] We have    and, for each $k$ even,      $m_k\le (3\max\{1, \ga\}/\sqrt{\al})^k \, k!!$. \\
\ite[d)] For each $k$, $\Var(N^{-1}\sum_{i=1}^N\lam_i^k)=O(N^{-1}).$\ent
\en{propo}

  The following theorem gives the limit local behavior of the $\lam_i$'s.
\beg{Th}[Poisson limit for bounded $N\bet$]\la{CVPOISSOn_smalllambda.6513}Suppose that $N\bet\lto 2\ga\ge 0$ as $N\to+\infty$ and that $\al\sim N\bet/2+1$. Fix $E \in \R$. As $N\to\infty$, the   point process $$\sum_{i=1}^N \del_{N(\lam_i-E)}$$ with $(\lam_1, \ld, \lam_N)\sim  P_{\al,\bet}^{(N)}$, converges in distribution to the law of a Poisson point process with intensity $\tta\,\ud x$ on $\R$, for \be\la{eq:theta(E)}\tta :=\f{(\ga+1)^{\ga+\f{1}{2}}}{\sqrt{2\pi}\Ga(\ga+1)}\exp\lf\{-\f{\ga+1}{2}E^2+2\ga\int\log|E-x|\ud\mu_\ga(x)\ri\}\,,\ee
with $\mu_\ga$ is in Proposition \re{Icnbty413}.
\en{Th}
   
   \beg{rmk}Note that the formula of $\tta$ given at \eqre{eq:theta(E)} should agree with the density of $\mu_\ga$ at $E$ as given in \cite{AllezGuionnetBouchaud,TrinhTomoyuki}, but we were not able to prove it so far.
   \en{rmk}
   
   \subsection{Case where $\bet\gg N^{-1}$}
  Using the fact that $N^2\bet\gg N$, one can easily adapt the proof of the following theorem from \cite{alice-greg-ofer}.
  \beg{Th}[LDP for $\bet\gg N^{-1}$]\la{mainth_165}Suppose that as $N\to\infty$, $\al=\al(N)$ and $\bet=\bet(N)$ are \st $\al\sim N\bet/2$,  $N\bet \gg {1}$ and $\bet$ is bounded. Then 
for $(\lam_1, \ld, \lam_N)$   distributed according to $P_{\al,\bet}^{(N)}$,    the sequence of random \pro measures $L_N:=\ff{N}\sum_{i=1}^N \del_{\lam_i}$ satisfies a LDP in the set of \pro measure on $\R$ endowed with the weak topology with speed $N^2\bet $ and good rate function $I$ defined by  \be\la{31111}I(\mu):=  \iint f(x,y)\mu\otd(\ud x \ud y)-\f{3}{8},\ee with $f: \R^2\to \R\cup\{+\infty\}$ is the function defined by   \be\la{31113}f(x,y)=\f{x^2+y^2}{8}-\ff{2}\log|x-y|.\ee Moreover, the unique minimum of $I$ is achieved at the semicircle law $$\si:=\ff{2\pi}\sqrt{4-x^2}\one_{|x|\le 2}\ud x$$  and we have  \be\la{lmct185}\lim_{N\to\infty} \ff{N^2\bet}\log Z_N(\al,\bet)=-\iint f(x,y)\si\otd(\ud x \ud y)=-\f{3}{8}.\ee
 \en{Th}
 
 In the case where $N\bet\gg 1$, as far as the local regime is concerned, we only have the following partial result, inspired from Johansson's work in \cite{Johansson}. Below, we explain how formally, it allows to prove the convergence of local statistics to the ones of a Poisson point process and to identify its density. 
 
\beg{Th}\la{PartialResultBetgg1}
Let $\bet=\bet(N)$ and $\al=\al(N) $ positive    \st \be\la{2231416h}\ff{N}\ll\bet\ll  \ff{\log N}\qquad;\qquad N\bet-2\al  \ll1.\ee Let $h~:~\R\to\R$ be a    bounded  function  having  9 continuous bounded derivatives   and $(\lam_1, \ld, \lam_N)$ be distributed according to $P_{\al,\bet}^{(N)}$. Then as $N\to\infty$,   we have \begin{equation}\label{Johlim}\log \int \mre^{\beta \sum_{j=1}^N h(\lam_j)}P^{(N)}_{\al,\bet}(\ud\lam_1, \ld, \ud\lam_N)-N\beta \int h (t) \ud\sigma(t)  \lto  \int h(t) \ud\nu (t)\end{equation} for \be\la{1731414h1a}\nu:=\ff{2}(\del_{-2}+\del_2)-\one_{|x|\le 2}\ff{\pi\sqrt{4-x^2}}\ud x\ee
\en{Th}

\beg{rmk}The measure $\nu$ of  \eqre{1731414h1a} is a classic correction to the semi-circle law (see e.g. \cite[Rem. 2.5]{Johansson} or, more recently, \cite{EnriquezMenard}).
\en{rmk}

Let us now explain how, on the formal level, Theorem \re{PartialResultBetgg1} gives, for any $E\in (-2,2)$, the convergence of the point process 
$$\sum_{i=1}^N \del_{N(\lam_i-E)}$$ to a Poisson point process with density $\tta \ud x$ on $\R$, for  \be  \la{203154DSO}\tta:= \f{\sqrt{4-x^2}}{2\pi},\ee 
The first thing one has to notice is that \be  \la{203154}\tta=\ff{2\pi}\exp\int \log|E-t|\ud\nu(t)\ee for $\nu$ being as in \eqre{1731414h1a} (the proof goes along the same lines as \cite[Lem. 2.7]{BAG3}).

To prove it rigorously, we would need to prove that for  \be\la{2702151bis}\tR^{(N)}_k(x_1, \ld, x_k)\ :=\ \int \mre^{\bet\sum_{i=1}^k\sum_{j=1}^{N-k}\log|E+\f{x_i}{N}-\lam_j| }P^{(N-k)}_{\al,\bet}(\ud\lam_{1}, \ld,\ud\lam_{N-k}),
\ee    for any $x_1, \ld, x_k\in \R$, 
 \be\la{eq:correlationfunctionsbis}  \f{Z_{N-k}(\al,\bet)}{Z_N(\al,\bet)} 
 \mre^{ - \f{\al}{2}\sum_{i}(E+\f{x_i}{N})^2}\tR^{(N)}_k(x_1, \ld, x_k)\;\lto\;\tta^k \ee and that we have an upper bound of the type of \eqre{861519h1}.
 
 First, it can be proved (see Section \re{2231420h560000}) that as $N\to\infty$,  for any fixed $k$,  \be\la{2231416h1}\f{Z_{N-k}(\al,\bet)}{Z_N(\al,\bet)}\sim \lf(\f{\mre^{N\bet/2}}{2\pi}\ri)^k.\ee
Moreover, there is a universal   positive  constant $M $ independent of $N$ \st uniformly on $N,k$,\be\la{2231416h1bis}\one_{1\le k\le N}\f{Z_{N-k}(\al,\bet)}{Z_N(\al,\bet)}\le  M^k\lf(\f{\mre^{N\bet/2}}{2\pi}\ri)^k.\ee
 
 Theorem \re{PartialResultBetgg1} can be rewritten as follows: for each fixed $h$ as in the theorem, 
\be\la{253149h} \E_{P^{(N)}_{\al,\bet}} \mre^{\beta \sum_{j=1}^N h(\lam_j)}=\exp\lf\{  N\beta \int h (t) \ud\sigma(t)+ \int h(t) \ud\nu (t)+\eps_N(h)\ri\},\ee with $\eps_N(h)\ll 1$.
 By \eqre{253149h}, cutting on the right thanks to Lemma \re{2231420h56} and making as if the function $h_N:\lam\mapsto \sum_{i=1}^k\log|E+\f{x_i}{N}-\lam| $ were $\Cc^9$ (and close enough to the function $h:\lam\mapsto k\log|E-\lam|$), we should have 
 \beq \tR^{(N)}_k(x_1, \ld, x_k)&=& \int \mre^{N\bet\sum_{i=1}^{N-k}h_N(\lam_j) }P^{(N-k)}_{\al,\bet}(\ud\lam_{1}, \ld,\ud\lam_{N-k})\\
 &=&
 \exp\lf\{  N\beta \int h_N (t) \ud\sigma(t)+ \int h_N(t) \ud\nu (t)+\eps_N(h_N)\ri\}\\
 &\approx&
 \exp\lf\{  N\beta \int h (t) \ud\sigma(t)+ \int h(t) \ud\nu (t)+\eps_N(h)\ri\}
 \eeq for $\nu$ as in \eqre{1731414h1a}.
 But    by \cite{BAG3} p. 529, we know that for any $E\in (-2,2)$, $$ \int k^{-1}h (t) \ud\sigma(t)=\int \log|E-t|\ud \si(t)=\f{E^2}{4}-\ff{2},$$ so that we should have 
\beqy\la{203151} \tR^{(N)}_k(x_1, \ld, x_k)&\approx&\exp\lf\{  N\beta \lf(\f{kE^2}{4}-\f{k}{2}\ri)+ \int h(t) \ud\nu (t)+\eps_N(h)\ri\}
 .\eeqy
Besides, by (\re{2231416h1}), we have \be\la{203152} \f{Z_{N-k}(\al,\bet)}{Z_N(\al,\bet)}\sim\lf(\f{\mre^{N\bet/2}}{2\pi}\ri)^k.\ee
Puting together \eqre{203151}, \eqre{203152} and the fact that $\al\sim\f{N\bet}{2}$, we should have 
 \beqy\nonumber\f{Z_{N-k}(\al,\bet)}{Z_N(\al,\bet)} 
 \mre^{ - \f{\al}{2}\sum_{i}(E+\f{x_i}{N})^2}\tR^{(N)}_k(x_1, \ld, x_k) &\approx & \lf(\f{\mre^{N\bet/2}}{2\pi}\ri)^k\mre^{-\f{kN\bet}{4}E^2}\mre^{  N\beta \lf(\f{kE^2}{4}-\f{k}{2}\ri)+  \int h(t) \ud\nu (t)}\\
 \la{203153} &=&\tta^k
 \eeqy
with $\tta$ as in \eqre{203154}.

\section{Proof of Proposition \re{Icnbty413}}
 Let $H$ be as in \eqre{matH6513}. We shall prove that for any $k$, $N^{-1}\E\Tr H^k$ tends to $m_k$, that the $m_k$'s satisfy Carleman's criterion   and that  d) holds. By Skorohod's representation theorem (see e.g. \cite[Th. 2.3.2]{Durrett2010})  and a diagonal extraction, it will imply a). Part c) will be clear from the proof.   Note first that if $Y$ is a   $\Ga(t)$-distributed variable, then for all $k\ge 0$, \be\la{momentsgamma20313} \E Y^{k}=\f{\Ga(t+k)}{\Ga(t)}.\ee

Let us fix $k\ge 0$, $u\in [0,1]$, let $i=i(N)$ be \st $i/N\lto u$ and compute $\E (H^{k})_{ii}$. We have $$(H^k)_{ii}=\sum_\eps H_{\eps(0), \eps(1)}\cdd H_{\eps(k-1),\eps(k)},$$ where the
sum is taken over paths $\eps:\{0, \ld,k\}\to \{1,Â \ld, n\}$ \st \bgt\ite $\eps(0)=\eps(k)=i$,\\
\ite for all $\ell=1, \ld, k$, $\eps(\ell)-\eps(\ell-1)=-1,0$ or $1$, in which case we say that $\ell$ belongs respectively to $\op{D}(\eps)$, $\op{F}(\eps)$ or $\op{U}(\eps)$.
\ent

Note first that for such a path $\eps$, 
\be\la{284131}\# \op{D}(\eps)=\# \op{U}(\eps), \qquad \# \op{D}(\eps)+\# \op{F}(\eps)+\# \op{U}(\eps)=k.\ee For any $j\in\{1, \ld, N\}$, we introduce \beq  \op{F}_j(\eps)&:=&\{\ell\in \{1, \ld, k\}\ste \eps(\ell-1)=\eps(\ell)=j\}\\
 \op{U}_j(\eps)&:=&\{\ell\in \op{U}(\eps)\ste  \eps(\ell)=j\}
\eeq
Then one can easily see, using \eqre{momentsgamma20313},  that  \beq \al^{k/2}\E H_{\eps(0), \eps(1)}\cdd H_{\eps(k-1),\eps(k)} &=&\prod_{j}(\E g^{\#\op{F}_j(\eps)} \E Y_j^{\#\op{U}_j(\eps)})\\ &=& \one_{\forall j, \;\#\op{F}_j(\eps)\trm{ is even}}\ti\\ &&\prod_{j}\lf( \#\op{F}_j(\eps)!!\ti \f{\Ga\lf(\f{(j-1)\bet}{2}+\#\op{U}_j(\eps)\ri)}{\Ga(\f{(j-1)\bet}{2})}\ri)\\
&\Ninf&  \one_{\forall j, \; \#\op{F}_j(\eps)\trm{ is even}}\ti\\ &&\prod_{j}\lf( \#\op{F}_j(\eps)!!\ti \ga u(\ga u+1)\cd (\ga u+  \#\op{U}_j(\eps)-1) \ri)
\eeq
By the dominated convergence theorem, it follows   that $N^{-1}\E\Tr H^k$ converges to  
\be\la{ictwyuuhb} m_k:= \sum_{\eps} (\ga+1)^{-k/2}\int_{u=0}^1\prod_{j}\lf\{ \#\op{F}_j(\eps)!!\ti \ga u(\ga u+1)\cd (\ga u+  \#\op{U}_j(\eps)-1)\ri\}\ud u,\ee where the sum runs over paths $\eps : \{0, \ld, k\}\to \z$ whose steps are in $\{-1, 0, 1\}$, \st $\eps(0)=\eps(k)=0$ and for all $k$, $\#\op{F}_j(\eps)$ is even.

Note that $\sum_j\#\op{F}_j(\eps)=\#\op{F}(\eps)$, whose parity is the one of $k$ by \eqre{284131}, so that when $k$ is odd, $\E (H^{k})_{ii}=0$. Using  \eqre{284131} again, we see that when $k$ is even, for any $\eps$, for any $u$, $$\prod_{j}\lf( \#\op{F}_j(\eps)!!\ti \ga u(\ga u+1)\cd (\ga u+  \#\op{U}_j(\eps)-1) \ri)\le \max\{1,\ga\}^k k!!,$$ so that the $m_k$'s satisfy Carleman's criterion. It follows that the $m_k$'s are the moments of a unique measure $\mu_\ga$ which depends continuously on $\ga$.  Besides, d) follows from the fact that 
  $\op{Cov}((H^k)_{ii},(H^k)_{jj})=0$ as soon as $|j-i|>2k$.

If $\ga=0$, then the only way for the term associated to $\eps$ in \eqre{ictwyuuhb}  to be non zero is that $k$ is even and $\eps$ is the constant path equal to $i$. This proves that $\mu_0=N(0,1)$. The reciprocal is obvious, as the fact that $\mu_\ga$ tends  to the semicircle law when $\ga\to\infty$ (using the formula of the moments of the semicircle law in terms of Dyck paths, as in \cite{alice-greg-ofer} or \cite{Ns06}).

To prove the first part of  b),  we use Lemma \re{cor324513} below. 
For any $a\in \R$ and $\eps\in (0,1]$, we have $$\mu((a-\eps, a+\eps))\le \liminf_{N\to\infty} P_{\al,\bet}^{(N)}(|\lam_1-a|< \eps)\le C\eps \mre^{\f{\ga a^2}{4} }.$$ The second part of b) is a direct consequence of Lemma \re{propodefmu75132} below.

\section{Proof of Theorem \re{CVPOISSOn_smalllambda.6513}}

\subsection{Correlation functions}\la{subsec:correlation_functions}
To prove the theorem, according to Proposition \re{prop:convPoisson}, we introduce the correlation functions of the point process $\sum_{i=1}^N\del_{N(\lam_i-E)}$, given by the formulas   \be\la{eq:correlationfunctions}R^{(N)}_k(x_1, \ld, x_k)\ :=\ \f{Z_{N-k}(\al,\bet)}{Z_N(\al,\bet)}\f{N^{-k}N!}{(N-k)!}\f{\Del(x_1, \ld, x_k)^\bet}{N^{\f{\bet k(k-1)}{2}}}\mre^{ - \f{\al}{2}\sum_{i}(E+\f{x_i}{N})^2}\tR^{(N)}_k(x_1, \ld, x_k)\ee
with \be\la{2702151}\tR^{(N)}_k(x_1, \ld, x_k)\ :=\ \int \mre^{\bet\sum_{i=1}^k\sum_{j=1}^{N-k}\log|E+\f{x_i}{N}-\lam_j| }P^{(N-k)}_{\al,\bet}(\ud\lam_{1}, \ld,\ud\lam_{N-k})
\ee

First of all, we know that \be\la{eq:ubfact} \one_{k\le N}\f{N^{-k}N!}{(N-k)!}\le 1\ee and that as $N\to\infty$, for each fixed $k$, \be\la{eq:cvfact} \f{N^{-k}N!}{(N-k)!}\Ninf 1.\ee
Besides, for any    $M>0$, for any $k\ge 1$ and any $x_1, \ld, x_k\in [-M,M]$, we have \be\la{eq:ubdelta} \one_{k\le N}\f{\Del(x_1, \ld, x_k)^\bet}{N^{\f{\bet k(k-1)}{2}}}\le \lf(\f{(2M)^\bet}{N^{\f{N\bet}{2}}}\ri)^k\ee
and  as soon as $\bet\ll(\log N)^{-1}$, for any fixed $k$ and any fixed $x_1, \ld x_k$, \be\la{eq:cvdelta} \f{\Del(x_1, \ld, x_k)^\bet}{N^{\f{\bet k(k-1)}{2}}}\Ninf 1.\ee
 
 \subsection{Partition functions}

We know, by \cite[Cor. 2.5.9]{alice-greg-ofer},   that\be\la{eq:ZNalphabeta} Z_N(\al, \bet) 
\ =\ \al^{-\lf(\f{N(N-1)}{4}\bet +\f{N}{2}\ri)}(2\pi)^{N/2} \prod_{j=1}^{N}\f{\Ga(1+j\bet/2)}{\Ga(1+\bet/2)}
\ee 
Hence 
\beq  \one_{k\le N} \f{Z_{N-k}(\al, \bet) }{ {Z}_N(\al, \bet) }&=& \al^{\f{k}{2}}\al^{\f{\bet}{4}(N(N-1)-(N-k)(N-k-1))}\prod_{j=N-k+1}^{N}(2\pi)^{-1/2}\f{ \Ga(1+\f{{\bet}}{2})}{\Ga(1+\f{j{\bet}}{2}) }\\
&=& \al^{\f{k}{2}}\al^{\f{\bet}{4}(2kN-k(k+1))}\prod_{j=N-k+1}^{N}(2\pi)^{-1/2}\f{ \Ga(1+\f{{\bet}}{2})}{\Ga(1+\f{j{\bet}}{2}) }
\eeq
By hypothesis,  $N\bet$ is bounded and so is $\al$. Let $C\ge 1$ be \st   $N\bet +\al\le C$. 
Then we have, uniformly in $k$, 
\be\la{eq:ubz} \one_{k\le N} \f{Z_{N-k}(\al, \bet) }{ {Z}_N(\al, \bet) }\le \lf(\f{C^{\f{1+C}{2}}}{\sqrt{2\pi}}\f{\max_{[1,  C]}\Ga}{\min_{[1,  C]}\Ga}\ri)^k\ee
Besides, as $N\bet\lto2\ga\ge 0$,   for each fixed $k$, \be\la{eq:cvz}\f{Z_{N-k}(\al, \bet) }{ {Z}_N(\al, \bet) }\lto \lf(\f{(\ga+1)^{\ga+\f{1}{2}}}{\sqrt{2\pi}\Ga(\ga+1)}\ri)^k.\ee

\subsection{Uniform upper-bound on the correlation functions}
 
 \beg{lem}\la{unifT_k6513}Let $\mc{K}$ be a compact subset of $\R$. There is a constant $C$ depending only on $\mc{K}$ and on the upper bounds on the sequences $N\bet,\al$ \st for all $k,N$ and all $x_1, \ld, x_k\in \mc{K}$, we have $$ \one_{k\le N}\tR_k^{(N)}(x_1, \ld, x_k)\le  C^k\mre^{\ff{4}\sum_{i=1}^k x_i^2}.$$
 \en{lem}
 
 \bpr  Note that by \eqre{eqpratique1} and \eqre{eqpratique2},  for each $i\in \{1, \ld k\}$ and $j\in \{1, \ld, N-k\}$, we have 
   $$|E+\f{x_i}{N}-\lam_j|^\bet\le (|E+\f{x_i}{N}|+|\lam_j|)^\bet\le 2^\bet \exp\{\bet \f{(E+\f{x_i}{N})^2+\lam_j^2}{8}\}.$$ Hence for $C$ a constant  (that might change from line to line) as in the statement of the lemma, 
   \beq \tR_k^{(N)}(x_1, \ld, x_k)&\le & 2^{kN\bet}\exp\{\f{N\bet}{8}\sum_{i=1}^k (E+\f{x_i}{N})^2\}\int\mre^{\f{k\bet}{8}\sum_{j=1}^{N-k} \lam_j^2}P_{\al,\bet}^{(N-k)}(\ud\lam_1, \ld, \ud\lam_{N-k})\\
   &\le&\f{Z_{N-k}(\al-k\bet/4,\bet)}{Z_{N-k}(\al ,\bet)}C^k\mre^{kCb^2}\\
   &\le& C^k\mre^{kCb^2}
   \eeq
 where we used    \eqre{eq:ZNalphabeta} and  the fact that for any $x\in [0,1/2]$, $(1-x)^{-1}\le 4^x$. \epr

Hence by \eqre{eq:ubfact}, \eqre{eq:ubdelta}, \eqre{eq:ubz} and the previous lemma, we have proved that b) of Proposition \re{prop:convPoisson} is satisfied. It remains to prove a) for $\theta$ given by \eqre{eq:theta(E)}.

\subsection{Preliminary estimates}

\beg{lem}[Bulk eigenvalues]\la{cor324513}There is a constant $C$ depending only on the upper bounds on the sequences  $\al$ and $N\bet$  \st for any $a\in \R$  and $\eps\in (0,1]$, $$P_{\al,\bet}^{(N)}(|\lam_1-a|\le \eps)\;\le\; C\eps   \mre^{\f{N\al\bet}{2(4\al-\bet)}a^2},$$
\en{lem}

\bpr 
We have \beq P_{\al,\bet}^{(N)}(|\lam_1-a|\le \eps)&= &\ff{Z_N(\al,\bet)}\int_{\lam_1\in[a\pm\eps]}\ud \lam_1\mre^{-\f{\al}{2}\lam_1^2} \int_{\mu\in\R^{N-1}}\Del(\mu)^\bet \prod_{j}|\lam_1-\mu_j|^\bet \mre^{-\f{\al}{2}\sum_j\mu_j^2}\ud \mu
\eeq
Note that if $|\lam_1-a|\le \eps$, the for any $j$, $|\lam_1-\mu_j| \le \eps+|\mu_j-a| .$ Moreover, for all $x\in \R$, $t>0$, we have  $|x|\le t^{-1/2}\mre^{\f{t x^2-1}{2}},$   hence \be\la{eqpratique1} |x|\le 2\mre^{\f{x^2}{16}}.\ee Using also the fact that \be\la{eqpratique2}(x+y)^2\le 2x^2+2y^2,\ee we get that if $|\lam_1-a|\le \eps$,$$|\lam_1-\mu_j|^\bet\le( \eps+|\mu_j-a|)^\bet\le 2^\bet \mre^{\bet\f{\eps^2+(\mu_j-a)^2}{8}}. $$
Hence 
\beq P_{\al,\bet}^{(N)}(|\lam_1-a|\le \eps)&\le &\f{2^{N\bet}\mre^{\f{N\bet\eps^2}{8}}}{Z_N(\al,\bet)}\int_{\lam_1\in[a\pm\eps]}\ud \lam_1\mre^{-\f{\al}{2}\lam_1^2} \int_{\mu\in\R^{N-1}}\Del(\mu)^\bet   \mre^{-\f{\al}{2}\sum_j(\mu_j^2-\bet\f{(\mu_j-a)^2}{4\al})}\ud \mu
\eeq
We have \beq \mu_j^2-\bet\f{(\mu_j-a)^2}{4\al}
%&=&\ff{4\al}\{(4\al-\bet)\mu^2+2\bet\mu a -\bet a^2\}\\
%&=&\f{4\al-\bet }{4\al}\{\mu^2+2\f{\bet a}{4\al-\bet}\mu  -\f{\bet a^2}{4\al-\bet}\}\\
%&=&\f{4\al-\bet }{4\al}\{(\mu+\f{\bet a}{4\al-\bet})^2   -(\f{\bet }{4\al-\bet}+\f{\bet^2 }{(4\al-\bet)^2})a^2\}\\
%&=&\f{4\al-\bet }{4\al}(\mu+\f{\bet a}{4\al-\bet})^2   -\f{4\al-\bet}{4\al}(\f{\bet }{4\al-\bet}+\f{\bet^2 }{(4\al-\bet)^2})a^2\\
%&=&\f{4\al-\bet }{4\al}(\mu+\f{\bet a}{4\al-\bet})^2   - (\f{\bet }{4\al}+\f{\bet^2 }{4\al(4\al-\bet)})a^2\\
&=&\f{4\al-\bet }{4\al}(\mu_j+\f{\bet a}{4\al-\bet})^2   - \f{\bet }{4\al-\bet}a^2
\eeq
Hence \beq P_{\al,\bet}^{(N)}(|\lam_1-a|\le \eps)&\le &\f{2^{N\bet}\mre^{\f{N\bet\eps^2}{8}} \mre^{\f{N\al\bet}{2(4\al-\bet)}a^2}}{Z_N(\al,\bet)}\int_{\lam_1\in[a\pm\eps]}\ud \lam_1\mre^{-\f{\al}{2}\lam_1^2} \int_{\mu\in\R^{N-1}}\Del(\mu)^\bet   \mre^{-\f{\al}{2}\sum_j\f{4\al-\bet}{4\al}(\mu_j+\f{\bet a}{4\al-\bet})^2}\ud \mu\\
&=& \eps 2^{N\bet+1} \mre^{\f{N\bet\eps^2}{8}} \mre^{\f{N\al\bet}{2(4\al-\bet)}a^2}  \f{Z_{N-1}(\al-\bet/4,\bet)}{Z_N(\al,\bet)}\\
&\le & C\eps \mre^{\f{N\al\bet}{2(4\al-\bet)}a^2},
\eeq
 we we used \eqre{eq:ZNalphabeta} to upper bound partition functions quotient. 
\epr

 \beg{lem}[Largest eigenvalues]\la{propodefmu75132} There is a constant $C$ depending only on the upper bounds on the sequences  $\al$ and $N\bet$  \st  for all $x>0$, \be\la{majgrandlambda_1}P_{\al,\bet}^{(N)}(|\lam_1|\ge x )  \;\le \; C\f{\mre^{- (4\al-\bet)\f{x^2}{8}}}{x} .\ee
\en{lem}

\bpr We have \beq P_{\al,\bet}^{(N)}( |\lam_1|\ge x )
&=& \ff{Z_N(\al,\bet)}\int_{|\lam_1|\ge x} \ud \lam_1 \mre^{-\f{\al}{2}\lam_1^2}  \int_{\mu\in \R^{N-1}} \Del(\mu)^\bet \prod_j|\lam_1-\mu_j|^\bet \mre^{-\f{\al}{2}\sum_j\mu_j^2}
\ud \mu\eeq
Note that by \eqre{eqpratique1} and \eqre{eqpratique2},    $$|\lam_1-\mu_j|^\bet\le (|\lam_1|+|\mu_j|)^\bet\le 2^\bet \mre^{\bet \f{\lam_1^2+\mu_j^2}{8}}.$$
 Hence \beq P_{\al,\bet}^{(N)}( |\lam_1|\ge x )
&\le& \f{2^{N\bet}}{Z_N(\al,\bet)}\int_{|\lam_1|\ge x} \ud \lam_1 \mre^{-\f{1}{2}(\al-\f{\bet}{4}) \lam_1^2 }  \int_{\mu\in \R^{N-1}} \Del(\mu)^\bet   \mre^{-\ff{2}(\al-\f{\bet}{4})\sum_j\mu_j^2}\ud \mu
\\
&=& \f{2^{N\bet}}{\sqrt{\al-\bet/4}}\f{ Z_{N-1}(\al-\bet/4,\bet)}{Z_N(\al,\bet)}\int_{|\lam_1|\ge x\sqrt{\al-\bet/4}}\mre^{-\f{\lam_1^2}{2}} \ud \lam_1   
\eeq Then we conclude using   \eqre{eq:ZNalphabeta} 
  and  the fact that  for   all $y>0$, \be\la{majGauss4513}\int_{y}^{+\infty} \mre^{-\f{t^2}{2}}\ud t\;\;\le\;\;  \f{\mre^{-\f{y^2}{2}}}{y} .\ee
\epr

%\subsection{Upper-bound for the tail of the empirical spectral law}\la{sec:Upper-bound for dev}

\beg{lem}[Tail of the empirical spectral law]There are  some   constants  $C,c $  depending only on the upper bounds on the sequences $N\bet, \al$ \st   for all $N $ and all $x>0$, $$P_{\al,\bet}^{(N)}(\f{\lam_1^2+\cdd+\lam_N^2}{N}\ge x)\;\le\; C\mre^{-c x}.$$
\en{lem}

\bpr We use again the tridiagonal matrix model of \eqre{matH6513} for $P_{\al,\bet}^{(N)}$  of the $\lam_i$'s. We know that  $\lam_1^2+\cd+\lam_N^2$ has the same law as $\Tr H^2$, with $H$ the matrix introduced as \eqre{matH6513}. Note that by the well known convolution relations between Gamma-distributed variables, $\ds\f{\al}{2}\Tr H^2$ has a $\Ga(\vfi)$-distribution for $\vfi:=N(1+(N-1)\bet)/2$. Hence  $$P_{\al,\bet}^{(N)}(\f{\lam_1^2+\cdd+\lam_N^2}{N}\ge x)\;\le \; \p(G\ge \f{\al Nx}{2}),$$
Then, one concludes using the concentration inequalities for Gamma variables (see \cite{BLM} p. 28-29) which say that for all $u\ge 0$, $$P(G\ge \vfi(1+u))\ \le \ \mre^{-\vfi(1+u-\sqrt{1+2u})}.$$ 
\epr

\beg{lem}\la{lem87131h57} For $C,c$ as in the previous lemma,   for any   $u,M,\tta$ \st $cM^2>2\tta\ge 0$, $$ 1\;\le\;\int\mre^{{\f{\tta}{N}}\sum_i\log(\f{|\lam_i-u|}{M}\vee 1)}P_{\al,\bet}^{(N)}(\ud\lam_1, \ld, \ud\lam_N)\;\le\; 1+\f{2C\tta \mre^{c u^2}}{cM^2-2\tta}.$$ 
\en{lem}

\bpr The integral above rewrites   \be\la{85131:06}  1+\tta\int_0^{+\infty}\!\!\!\!\mre^{\tta x}P_{\al,\bet}^{(N)}(\ff{N}\sum_i\log(\f{|\lam_i-u|}{M}\vee 1)\ge x)\ud x.\ee
Now, note that as $\log(|\lam-u|\vee 1)\le (\lam-u)^2$,  \beq P_{\al,\bet}^{(N)}(\ff{N}\sum_i\log(\f{|\lam_i-u|}{M}\vee 1)\;\ge\; x)&\le &P_{\al,\bet}^{(N)}(\f{(\lam_1-u)^2+\cdd+(\lam_N-u)^2}{N}\ge M^2x)\\
&\le &P_{\al,\bet}^{(N)}(\f{2\lam_1^2+2u^2+\cdd+2\lam_N^2+2u^2}{N}\ge M^2x)\\ 
&\le &P_{\al,\bet}^{(N)}(\f{\lam_1^2 +\cdd+\lam_N^2}{N}\ge \f{M^2x}{2}-u^2).\eeq Then one concludes using  \eqre{85131:06} and the previous lemma.
\epr

\subsection{Convergence of the correlation functions}Let us now prove a) of Proposition \re{prop:convPoisson}  for $\theta$ given by \eqre{eq:theta(E)}.
Note first that by b) of Proposition \re{Icnbty413}, we know that    \be\la{eq:logL1muga}\int \log|E-x|\ud\mu_\ga(x)\;<\;\infty  .\ee 

Besides, by \eqre{eq:correlationfunctions}, \eqre{eq:cvfact}, \eqre{eq:cvdelta} and \eqre{eq:cvz}, it suffices to prove that for each $k$ and each $x_1, \ld, x_k\inÊ\R$, the quantity $\tR^{(N)}_k(x_1, \ld, x_k)$ defined at \eqre{2702151} satisfies, as $N\to\infty$, \be\la{eq:limtRk}\tR^{(N)}_k(x_1, \ld, x_k)\;\Ninf \; \exp\lf\{ 2\ga k\int\log|E-x|\ud\mu_\ga(x)\ri\}\,.\ee

 \subsubsection{Upper-bound}\la{lem87131h572}
 Let us prove that for any fixed $k$ and  $x_1, \ld, x_k$, \be\la{751300:33}\limsup_{N\to\infty} \tR^{(N)}_k(x_1, \ld, x_k)\;\le\; \exp\lf\{ 2\ga k\int\log|E-x|\ud\mu_\ga(x)\ri\}\,.\ee
For $\eps>0$, set
$$ \tR^{(N,\eps)}_k(x_1, \ld, x_k)\ :=\ \int \mre^{\bet\sum_{i=1}^k\sum_{j=1}^{N-k}\log(|E+\f{x_i}{N}-\lam_j|\vee\eps) }P^{(N-k)}_{\al,\bet}(\ud\lam_{1}, \ld,\ud\lam_{N-k}).$$ 
As $\ds \tR^{(N)}_k\le  \tR^{(N,\eps)}_k$ and, by \eqre{eq:logL1muga},  $\ds\inf_{\eps>0}\int\log(|x-E|\vee \eps)\ud \mu_\ga(x)=\int\log|x-E| \ud \mu_\ga(x)$, it suffices to prove that for any $\eps>0$ small enough, we have \be\la{751300:33bis}\limsup_{N\to\infty} \tR^{(N,\eps)}_k(x_1, \ld, x_k)\;\le\; \exp\lf\{ 2\ga k\int\log(|E-x|\vee \eps)\ud\mu_\ga(x)\ri\}\,.\ee

  Note now   that for any $M>0$ large enough, as $$\log(|E+x_i/{N}-x|\vee \eps)= \log\{(|E+x_i/{N}-x|\vee \eps)\wedge M\}+\log\lf(\f{|E+x_i/{N}-x|}{M}\vee 1\ri)$$ the function $ \tR^{(N,\eps)}_k(x_1, \ld, x_k)$ rewrites     \begin{eqnarray*}&&
 \tR^{(N,\eps)}_k(x_1, \ld, x_k)\\&=& \int  \mre^{\bet \sum_{i=1}^k \sum_{j=1}^{N-k} \log\{(|E+\f{x_i}{N}-\lam_j|\vee \eps)\wedge M\}} \mre^{\bet \sum_{i=1}^k \sum_{j=1}^{N-k} \log\lf(\f{|E+\f{x_i}{N}-\lam_j|}{M}\vee 1\ri)} P^{(N-k)}_{\al,\bet}(\ud\lam_{1}, \ld,\ud\lam_{N-k})\\
 &=& \int  \mre^{\bet \sum_{i=1}^k \sum_{j=1}^{N-k} \log\{(|E+\f{x_i}{N}-\lam_j|\vee \eps)\wedge M\}}  P^{(N-k)}_{\al,\bet}(\ud\lam_{1}, \ld,\ud\lam_{N-k})+\\
 && \int  \underbrace{\mre^{\bet \sum_{i=1}^k \sum_{j=1}^{N-k} \log\{(|E+\f{x_i}{N}-\lam_j|\vee \eps)\wedge M\}}}_{:=X} \underbrace{\lf(\mre^{\bet \sum_{i=1}^k \sum_{j=1}^{N-k} \log\lf(\f{|E+\f{x_i}{N}-\lam_j|}{M}\vee 1\ri)}-1\ri)}_{:=Y}P^{(N-k)}_{\al,\bet}(\ud\lam_{1}, \ld,\ud\lam_{N-k})
  \eeq
  
  By Proposition \re{Icnbty413}, we know that under the law $P^{(N-k)}_{\al,\bet}$,   the $L^\infty$-bounded  sequence of random variables   $\ds\bet \sum_{i=1}^k\sum_{j=1}^{N-k} \log\{(|E+\f{x_i}{N}-\lam_j|\vee \eps)\wedge M\}$ converges in probability, as $N\to\infty$,  to $$2\ga k\int\log\{(|E- x|\vee \eps)\wedge M\}\ud\mu_\ga(x)$$ (one gets rid of the $\f{x_i}{N}$'s by noticing, for example, that the convolution of \pro measures is continuous with respect to the weak topology and that $\del_{\lam_i/N}$ converges to $\del_0$).

  Note that by choosing $M$ large enough, one can make $\ds \int\log\{(|E- x|\vee \eps)\wedge M\}\ud\mu_\ga(x)$ as close as we want from $\int\log(|E-x|\vee \eps)\ud\mu_\ga(x)$. Moreover, one can easily adapt the proof of Lemma \re{unifT_k6513} to see that  $\int X^2P^{(N-k)}_{\al,\bet}(\ud\lam_{1}, \ld,\ud\lam_{N-k})$ is bounded by a constant independent of $M$, hence by Cauchy-Schwartz,  to prove  \eqre{751300:33bis}, it suffices to prove that $\int Y^2P^{(N-k)}_{\al,\bet}(\ud\lam_{1}, \ld,\ud\lam_{N-k})$ can be made as small as we want if $M$ is large enough.
 Note that for $$Y_i:= \mre^{\bet  \sum_{j=1}^{N-k} \log\lf(\f{|E+\f{x_i}{N}-\lam_j|}{M}\vee 1\ri)},$$ we have $Y=Y_1\cd Y_k-1$. Thus by the H\"older inequality, it is enough to prove that for $M$ large enough, earh $Y_i$ can have its $k$-th and $2k$-th moment as close as we want from $1$, which is a direct consequence of Lemma \re{lem87131h57}.

 \subsubsection{Lower bound}
 To obtain the analogous lower bound \be\la{751300:33lower}\liminf_{N\to\infty}  \tR^{(N)}_k(x_1, \ld, x_k)\;\ge\; \exp\lf\{ 2\ga k\int\log|E-x|\ud\mu_\ga(x)\ri\}\,,\ee
 we observe that first by Jensen's inequality and then by exchangeability,
\beq \log \tR^{(N)}_k(x_1, \ld, x_k)&\ge&
\int \bet\sum_{i=1}^k\sum_{j=1}^{N-k} \log |E+\frac{x_i}{N}-\lam_j|P^{(N-k)}_{\al,\bet}(\ud\lam_{1}, \ld,\ud\lam_{N-k})\\
&=& (N-k)\bet\int\sum_{i=1}^k  \log |E+\frac{x_i}{N}-\lam_1|P^{(N-k)}_{\al,\bet}(\ud\lam_{1}, \ld,\ud\lam_{N-k})
\eeq 
 Hence as $(N-k)\bet\lto 2\ga$ and the triplet $(N,\bet(N), \al(N))$ satisfies the same hypotheses as $(N-k,\bet(N), \al(N))$,  it suffices to prove that for    any fixed $x$,   we have  $$\liminf_{N\to\infty} \int  \log |E+\frac{x}{N}-\lam_1|P^{(N)}_{\al,\bet}(\ud\lam_{1}, \ld,\ud\lam_{N}) \;\ge\; \int \log|E-t|\ud\mu_\ga(t).
 $$
 As, by exchangeability, $\mu_\ga$ is also the weak limit of the distribution of $\lam_1$ under $P_{\al,\bet}^{(N)}$, 
 we know that for any $\eps>0$, $$\liminf_{N\to\infty} \int\log(|E+\f{x}{N}-\lam_1|\vee \eps)P^{(N)}_{\al,\bet}(\ud\lam_{1}, \ld,\ud\lam_{N}) \;\ge\; \int \log|E-t|\ud\mu_\ga(t)
 $$(and one can get rid  of   $\f{x}{N}$ for the same reason as in Section \re{lem87131h572} above).
 Hence it suffices that for $\eps$ small enough, $$\limsup_{N\to\infty}  \lf|\int\lf(\log(|E+\f{x}{N}-\lam_1|\vee \eps)  - \log|E+\f{x}{N}-\lam_1|\ri) P^{(N)}_{\al,\bet}(\ud\lam_{1}, \ld,\ud\lam_{N})\ri|$$ can be made as small as desired. But for any random variable $X>0$, $$\E[\log(X\vee \eps)-\log(X)]=\E[(\log\eps-\log(X))\one_{X\le \eps}]=\int_0^\eps\f{\p(X\le t)}{t}\ud t.$$ Here, by Lemma \re{cor324513},   there is a constant $C$ \st  $$P_{\al,\bet}^{(N)}(|E+\f{x}{N}-\lam_1|\le t)\;\le \; Ct,$$ which allows to get the desired bound.

 \section{Case where $N^{-1}\ll \beta \ll 1/\log(N)$}

\subsection{Partition functions: proofs of \eqre{2231416h1} and \eqre{2231416h1bis}}\la{2231420h560000}

    It follows from \eqre{eq:ZNalphabeta} that for all $N$, for all $1\le k\le N$, 
 \beq\f{Z_{N-k}(\al,\bet)}{Z_N(\al,\bet)}&=& \al^{\f{k}{2}}\al^{\f{\bet}{4}(N(N-1)-(N-k)(N-k-1))}(2\pi)^{-k/2}\Ga(1+\f{\bet}{2})^k\prod_{\ell=0}^{k-1}\ff{\Ga(1+(N-\ell)\bet/2)}\\
 &=& \al^{\f{k}{2}}\al^{\f{\bet}{4}( 2kN-k^2+k)}(2\pi)^{-k/2}\Ga(1+\f{\bet}{2})^k\prod_{\ell=0}^{k-1}\ff{\Ga(1+(N-\ell)\bet/2)} \eeq
 But by the Stirling formula, for $z\to+\infty,$ \be\la{fG3314}
\Ga(1+z)\sim \sqrt{2\pi z}\exp\lf\{z\log(z)	-z   \ri\}.\ee
Note that our hypothesis on $\bet$ implies that  for any fixed $\ell=0, \ld, k-1$, we have $$(N-\ell)\bet/2\log((N-\ell)\bet/2 )	-(N-\ell)\bet/2   =N\bet/2\log(N\bet/2)	-N\bet/2 +o(1),$$ so that, as \eqre{2231416h} implies that $N\bet\log\f{2\al}{N\bet}  \ll1$ and $\bet\log\al\ll 1$,   
\beq\f{Z_{N-k}(\al,\bet)}{Z_N(\al,\bet)}&\sim&\lf(\f{\mre^{N\bet/2}}{2\pi}\ri)^k\ti \lf(\exp\lf\{\f{N\bet}{2}\lf[\log\al -\log(N\bet/2)\ri]\ri\}\ri)^k
\eeq
By \eqre{2231416h}, we get \eqre{2231416h1}. The upper bound \eqre{2231416h1bis} comes in the same way, noticing that the error in \eqre{fG3314} is uniformly bounded on $z\ge 0$.

\subsection{Tail estimate}
\beg{lem}\la{2231420h56}Let $\bet=\bet(N)$ and $\al=\al(N) $ be   satisfying  \eqre{2231416h}. Then there is a constant $C $ depending only on the sequences $\al$ and $\bet$  \st for all $N$, for all $x>0$, \be\la{223142}P_{\al,\bet}^{(N)}(|\lam_1|\ge x )  \;\le \; C\mre^{\lf(- \f{x^2}{2} +C\ri)N\bet/2} .\ee
\en{lem}

\bpr We have \beq P^{(N)}_{\al,\bet}( |\lam_1|\ge x )
&=& \ff{Z_N(\al, \bet)}\int_{|\lam_1|\ge x} \ud \lam_1\mre^{-\f{\al}{2}\lam_1^2}  \int_{\mu\in \R^{N-1}} \Del(\mu)^\bet \prod_j|\lam_1-\mu_j|^\bet\mre^{-\f{\al}{2}\sum_j\mu_j^2}
\ud \mu\eeq
Note that by \eqre{eqpratique1} and \eqre{eqpratique2},    $$|\lam_1-\mu_j|^\bet\le (|\lam_1|+|\mu_j|)^\bet\le 2^\bet\mre^{\bet \f{\lam_1^2+\mu_j^2}{8}}.$$
 Hence \beq P^{(N)}_{\al,\bet}( |\lam_1|\ge x )
&\le& \f{2^{N\bet}}{Z_N(\al,\bet)}\int_{|\lam_1|\ge x} \ud \lam_1\mre^{-\f{\al-\bet/4}{2}  \lam_1^2 }  \int_{\mu\in \R^{N-1}} \Del(\mu)^\bet  \mre^{-\f{\al-\bet/4}{2} \sum_j\mu_j^2}\ud \mu
\\
&=& \f{2^{N\bet} }{Z_N(\al, \bet)\sqrt{\al-\bet/4}}\int_{|\lam_1|\ge x\sqrt{\al-\bet/4}}\mre^{-\f{\lam_1^2}{2}} \ud \lam_1\ti Z_{N-1}(\al-\f{\bet}{4},\bet)\\
 &\le &\f{2^{N\bet+1}\mre^{-\f{x^2(\al-\bet/4)}{2}}}{x(\al-\bet/4)}\f{Z_{N-1}(\al-\f{\bet}{4},\bet)}{Z_N(\al, \bet)}\\
 &=& \f{2^{N\bet+1}\mre^{-\f{x^2(\al-\bet/4)}{2}}}{x(\al-\bet/4)}(1-\bet/(4\al))^{-\lf(\f{(N-1)(N-2)}{4}\bet+\f{N-1}{2}\ri)}\f{Z_{N-1}(\al ,\bet)}{Z_N(\al, \bet)}
\eeq
Let us now use \eqre{2231416h1bis} and for example the fact that $(1-y)^{-1}\le\mre^{2y}$ when $y\in [0,1/2]$. We get 
 \beq P^{(N)}_{\al,\bet}( |\lam_1|\ge x )
&\le &\f{M}{2\pi}\ff{x\al} \exp\{ -\f{x^2}{2}(\al-\bet/4)+ N^2\bet^2/(8\al) +N\bet/(2\al) +N\bet/2\},
\eeq
which allows to conclude, as we already noticed that  \eqre{2231416h} implies that $2\al\sim N\bet\gg 1$.
\epr

 \subsection{Proof of Theorem \re{PartialResultBetgg1}}

We first  define the \pro measure on $\R^N$  \be\la{631414hdens}P^{(N,h)}_{\al,\bet}(\ud x):=\underbrace{\ff{Z_N^h(\al, \bet)}\Del(x)^\bet\exp\{-\f{\al}{2}\sum_{i=1}^N x_i^2+\bet\sum_{i=1}^N h(x_i)\}}_{\ds :=\rho^{(N,h)}_{\al,\bet}(x_1, \ld, x_N)}\ud x_1\cd\ud x_N,\ee where    \be\la{631414h}Z_N^h(\al, \bet)=\int \Del(x)^\bet\mre^{-\f{\al}{2}\sum_{i=1}^N x_i^2+\bet\sum_{i=1}^N h(x_i)} \ud x_{ 1}\cd \ud x_N\ee is the normalisation constant.
Let, for $i=1,2$, $$u_N^{i,h}(x_1, \ld, x_i):=  \int_{x_{i+1}, \ld, x_N}\rho^{(N,h)}_{\al,\bet}(x_1, \ld, x_N)\ud x_{i+1}\cd \ud x_N$$ be the $i$-th correlation  function of $\rho^{(N,h)}_{\al,\bet}$.  
\beg{lem}
 Let $\psi:\R\to\C$ be a $\Cc^1$ function on $\R$ \st the real and imaginary parts of $\psi'$ are bounded below. Then we have \beqy\nonumber \bet\f{N-1}{2}\iint  \f{\psi(t)-\psi(s)}{t-s}u_N^{2,h}(s,t)\ud s\ud t+ \int_{\R}(\bet h'(t)-\al t)\psi(t)u_N^{1,h}(t)\ud t&&\\ +\int_{\R} \psi'(t)u_N^{1,h}(t)\ud t&=&0\la{631414hbis2}.\eeqy
 \en{lem}
 
 \bpr As \eqre{631414hbis2} is linear in $\psi$, one can suppose   $\psi$ to be real-valued. Then for $  \tta\ge 0$ small enough,   the function $y+\tta\psi(y)$ is an homeomorphism on $\R$, hence one can make the change of variable $x_i=y_i+\tta\psi(y_i)$ in \eqre{631414h}. We get 
 \beqy\la{631414hbis}Z_N^h&=&\int_{y\in \R^N} \prod_{1\le i<j\le N}|y_j-y_i+\tta(\psi(y_j)-\psi(y_i))|^\bet \\ \nober&&e^{-\f{\al}{2}\sum_{i=1}^N (y_i^2+2\tta y_i\psi(y_i)+\tta^2\psi(y_i)^2)+\bet\sum_{i=1}^N h(y_i+\tta\psi(y_i))} \prod_{i=1}^N(1+\tta\psi'(y_i))\ud y_{ 1}\cd \ud y_N\eeqy
Let us compute the derivative, with respect to $\tta$, at $\tta=0$, of the RHT of \eqre{631414hbis}. We have $$\pa_{\tta, \tta=0}|y_j-y_i+\tta(\psi(y_j)-\psi(y_i))|^\bet=\bet\f{\psi(y_j)-\psi(y_i)}{y_j-y_i}|y_j-y_i|^\bet,$$ we also have
$$\pa_{\tta, \tta=0}\mre^{-\f{\al}{2}\sum_{i=1}^N (y_i^2+2\tta y_i\psi(y_i)+\tta^2\psi(y_i)^2)}= - \al\sum_{i=1}^Ny_i\psi(y_i),$$ $$\pa_{\tta, \tta=0}\mre^{\bet\sum_{i=1}^N h(y_i+\tta\psi(y_i))}=  \bet\sum_{i=1}^Nh'(y_i)\psi(y_i),$$ and $$\pa_{\tta, \tta=0}\prod_{i=1}^N(1+\tta\psi'(y_i))=  \sum_{i=1}^N \psi'(y_i).$$ Hence \beq &&\f{\pa_{\tta, \tta=0}(\trm{RHT of \eqre{631414hbis}})}{Z_N^h}\\&=&\bet\sum_{1\le i<j\le N}\int_{y\in \R^N}\f{\psi(y_j)-\psi(y_i)}{y_j-y_i}\rho_N^h(y_1, \ld, y_N)\ud y_1\cd \ud y_N\\
&&-\al\sum_{i=1}^N\int_{y\in \R^N}y_i\psi(y_i) \rho_N^h(y_1, \ld, y_N)\ud y_1\cd \ud y_N\\
&&+ \bet\sum_{i=1}^N\int_{y\in \R^N}h'(y_i)\psi(y_i) \rho_N^h(y_1, \ld, y_N)\ud y_1\cd \ud y_N\\
&&+\sum_{i=1}^N\int_{y\in \R^N} \psi'(y_i) \rho_N^h(y_1, \ld, y_N)\ud y_1\cd \ud y_N
\eeq
We get exactly \eqre{631414hbis2}.\epr

Now, we define, for $z$ \st $\Im z>0$, $$H_N(z):= \int\f{h'(t)}{z-t}u^{1,h}_N(t)\ud t\qquad ;\qquad U_N(z):= \int\ff{z-t}u^{1,h}_N(t)\ud t$$
 $$k_N(s,t):= Nu^{1,h}_N(s)u^{1,h}_N(t)-(N-1)u^{2,h}_N(s,t)$$ and 
 \beqy\la{661418h} K_N(z)&:=&N\iint\lf(\ff{(z-t)(z-s)}-\ff{2(z-t)^2}-\ff{2(z-s)^2}\ri)k_N(s,t)\ud t\ud s
\\ \nonumber
&=&
-\ff{2}\iint \f{(t-s)^2}{(z-t)^2(z-s)^2}k_N(s,t)\ud s\ud t.\eeqy
We also introduce
\be\la{f117314}U(z): = \ff{2}(z-\sqrt{z^2-4})\qquad\qquad (\Im(z)>0),\ee  where when $\Re(z)\ge 0$   (resp.  $\le 0$) , 
$\sqrt{z^2-4}$ is computed with the determination of the square root on $\C\bck(-\infty, 0)$ (resp. on $\C\bck(0,+\infty)$) with positive values on $[0,+\infty)$ (resp.  \st $\sqrt{-1}=\ii$). It is well known that $U$ is the Stieltjes transform of the semicircle law   $\si$.

\beg{lem}On the upper half-plane, we have  \be\la{1431417h}N\bet(U_N-U)(2U-\f{2\al}{N\bet}z+U_N-U)=\Del_N:=\bet N^{-1}K_N-2\bet H_N+(2-\bet)U_N'+(2\al-N\bet)U^2.\ee
\en{lem}

\bpr
We shall apply the previous lemma with  $\ds\psi(t)=\ff{z-t}$. Note that $\ds \f{\psi(t)-\psi(s)}{t-s}=\ff{(z-t)(z-s)},$ so that we have :
\beq &&N(N-1)\iint  \f{\psi(t)-\psi(s)}{t-s}u_N^{2,h}(s,t)\ud s\ud t\\
&=&\iint \ff{(z-t)(z-s)}N(N-1)u_N^{2,h}(s,t)\ud s\ud t\\
&=& \iint \ff{(z-t)(z-s)}(-Nk_N (s,t)+N^2u^{1,h}_N(s)u^{1,h}_N(t))\ud s\ud t\\
&=&
N^2U_N(z)^2-N\iint \ff{(z-t)(z-s)} k_N (s,t)\ud s\ud t\\
&=&N^2U_N(z)^2-N\iint \lf(\ff{(z-t)(z-s)}-\ff{2(z-t)^2}-\ff{2(z-s)^2}\ri) k_N (s,t)\ud s\ud t\\ &&-N\iint \ff{(z-t)^2} k_N (s,t)\ud s\ud t\\
&=&N^2U_N(z)^2-K_N(z)+nU_N'(z)\eeq
(where we use the fact that for any function $f(t)$, $\ds\iint f(t)k_N(s,t)\ud s\ud t=\int f(t)u_N^{1,h}(t)\ud t$). We also have  $$\int-t\psi(t)u_N^{1,h}(t)\ud t=1-zU_N(z)\quad;\quad \int\psi'(t)u_N^{1,h}(t)\ud t=-U_N'(z),$$ so $2N^{-2}\bet^{-1}\ti $\eqre{631414hbis2} rewrites 
$$U_N(z)^2-N^{-2}K_N(z)+N^{-1}U_N'(z)+2N^{-1}H_N(z)+\f{2\al}{N\bet}(1-zU_N(z))-\f{2 }{N\bet}U_N'(z)=0 $$
\ie  \be \la{631416h}U_N(z)^2-\f{2\al}{N\bet}zU_N(z) +\f{2\al}{N\bet}=N^{-2}K_N(z) -2N^{-1}H_N(z)+(\f{2 }{N\bet}-\ff{N})U_N'(z).\ee 
One gets \eqre{1431417h}, using the well known equation $U(z)^2-zU(z)+1=0$  (see \cite[Eq. (2.4.6)]{alice-greg-ofer}).\epr
 
A key step in the proof of the theorem will be  to prove that as $N\to\infty$, 
\be\la{keystep20314}\beta K_N(z) \ll N.\ee

We shall now  prove \eqre{keystep20314}. Let  $(y_1, \ld, y_N)$ be a random vector with distribution $P^{(N,h)}_{\al,\bet}$ and  for $g\in \Cc_b(\R, \C)$, define the random variable $$\hmu_N^h(g):=\sum_{i=1}^N g(y_i), $$ with variance $\ds   \Si_N(g):= \E [|\hmu_N^h(g)-\E [\hmu_N^h(g)]|^2]$.
As  \be\la{661418h1} K_N(z)=\E [(\hmu_N^h(g)-\E [\hmu_N^h(g)])^2]\ee for $g(t)=(z-t)^{-1}$, we have  \be\la{1431418h2}|K_N(z)|\le \Si_N(\ff{z-t}).\ee

  Note that as    \beq  \E [\hmu_N^h(g)] &=&N\int u_N^{1,h}(t)g(t)\ud t\\ \E [|\hmu_N^h(g)|^2]&=&\sum_{i,j=1}^N\E[g(y_i)\ovl{g(y_j)}]\\
  &=&N(N-1)\iint g(s)\ovl{g(t)}u^{2,h}_N(s,t)\ud s\ud t+N\int u_N^{1,h}(t)|g(t)|^2\ud t,
  \eeq we deduce  \beq \Si_N(g)&=& N(N-1)\iint g(s)\ovl{g(t)}(u^{2,h}_N(s,t)-u_N^{1,h}(s)u_N^{1,h}(t))\ud s\ud t\\
  &&+N\lf(\int u_N^{1,h}(t)|g(t)|^2\ud t-|\int u_N^{1,h}(t)g(t)\ud t|^2\ri)\eeq

\beg{lem}There is $L>0$ and $c>0$ \st for any fixed function $g$,  
\beq \Si_N(g)&\le & N(N-1)\iint_{[-L,L]^2} g(s)\ovl{g(t)}(u^{2,h}_N(s,t)-u_N^{1,h}(s)u_N^{1,h}(t))\ud s\ud t\\
  &&+N\lf(\int_{-L}^L u_N^{1,h}(t)|g(t)|^2\ud t-|\int_{-L}^L u_N^{1,h}(t)g(t)\ud t|^2\ri)+4N^2\|g\|_\infty\mre^{-cN\beta}\eeq
%$$\mathbb{P}\left (\exists i, |\lambda_i| >L/2\right )\leq\mre^{-cN\beta}.$$
\en{lem}

\bpr Using the fact that $|\sum_{i=1}^N h(x_i)|\le N\|h\|_\infty$, we see that the \pro measure $P^{(N,h)}_{\al,\bet}$ defined at \eqre{631414hdens} and its normalization constant can be controlled thanks to the \pro measure $P^{(N)}_{\al,\bet}$ and its normalization constant: for any Borel set $A\subset \R^N$, we have   $$P^{(N,h)}_{\al,\bet}(A)\le\mre^{2N\bet \|h\|_\infty}P^{(N)}_{\al,\bet}(A).$$ It follows that up to a change of the constant $C$, Lemma \re{2231420h56} is also true for $P^{(N,h)}_{\al,\bet}$, which allows to conclude.
\epr

This lemma allows to reduce the problem to a compact set, and after rescaling, one can turn the  compact set in question to  $[-1/2,1/2]$ : we deduce, as in \cite{Johansson}, that for $w:=z/L$ and  
$$\rho_N^2 (t,s):= u_N^{2,h} (Lt, Ls),\qquad  \rho_N^1 (t):=  u_N^{1,h}(Lt).$$
we have 
\be\la{1431418h}\Sigma_N(\frac{1}{z-t})\leq 
N(N-1) \int_{-1/2}^{1/2}\frac{1}{w-t}G_N(x\mapsto\frac{1}{\overline{w}-x})(t) \ud t+4N^2e^{-cN\beta}(\Im z)^{-2},\ee
where $G_N$ is the operator on $L^2 ([-1/2,1/2], \ud x)$ defined by $$
G_N(f) (t)=\int f(s) (\rho_N^2 (t,s)-\rho_N^1(t)\rho_N^1(s))\ud s+\frac{\rho_N(t)}{N-1}\left( \frac{1}{L}f(t)-\int_{-1/2}^{1/2}\rho_N (s) f(s) \ud s\right).
$$

Thus to prove the estimate of interest \eqre{keystep20314},  we have to upper bound:
\begin{equation}\label{rest}N(N-1) \int_{-1/2}^{1/2}\frac{1}{w-t}G_N(x\mapsto \frac{1}{\overline{w}-x})(t)\ud t.\end{equation}
Following \cite{Johansson}, we introduce the integral operator $P_w$ on $L^2 ([-1/2,1/2], \ud x)$
with kernel $P_w(t,s)=\frac{1}{(w-t)(\overline{w}-s)}.$ Then $(s,t)\longmapsto \frac{1}{w-t}G_N(x\mapsto \frac{1}{\overline{w}-x})(s)$ is an integral operator satisfying the hypothesis of Theorem 2.12 of \cite{SimonTrace}. This   trace class operator is nothing but $P_wG_N$, thus   by this theorem, we have 
  \be\la{1431418h1} \int_{-1/2}^{1/2}\frac{1}{w-t}G_N(\frac{1}{\overline{w}-s})(t)\ud t=   \Tr(P_w G_N).\ee
 We will not here recall all the arguments used in \cite{Johansson} Lemma 3.12 and Proposition 3.9 
to estimate this trace. The proof transfers to our setting using minor modifications (essentially replacing $N$  by $N\beta $).
Note that the important feature of $h$ is that is Lipschitz on compact sets in our case.
Thus we simply state the final estimate we will use in this article, namely the following lemma : 
\beg{lem} We have 
$\Tr(P_w G_N)\leq C N^{-1}\log(N)$ for some constant $C$. 
\end{lem}
It follows, by \eqre{1431418h2},  \eqre{1431418h} and \eqre{1431418h1}, that 
$$|\bet K_N(z)|\le\bet \Sigma_N(\frac{1}{z-t})\leq  C N\bet\log(N)+ +4N^2\bet\mre^{-cN\beta}(\Im z)^{-2}.$$ 
This is of course $\ll N$, so the estimate of interest \eqre{keystep20314} is proved.

As $\bet\log(N)\ll 1$ and $N\bet-2\al\ll 1$, by \eqre{1431417h}, we deduce   that, uniformly on compact subsets of $\C^+$, one has that 
\begin{equation}\la{1731418h}
N\bet (U_N(z)-U(z))\lto \frac{2U'(z)}{ 2U(z)-z} =\f{z}{z^2-4}-\ff{\sqrt{z^2-4}}.\ee One recognizes easily that the RHT of \eqre{1731418h} is the Stieltjes transform of the null mass signed measure  $\nu $ of \eqre{1731414h1a}  (to do that, use the fact that $U(z)$, given by \eqre{f117314},  is the Stieltjes transform of the semi-circle law and then use an integration by parts).

The rest of the proof of the theorem is   an easy adaptation of  the proof of Theorem 2.4 in \cite{Johansson} (p. 169-172).
 The main idea is to define $$F(\lam):=\log \E\mre^{\beta \lam \sum_{j=1}^N h(x_j)}-N\beta \lam\int h (t) \ud\sigma(t),$$   to notice that 
 $$\pa_\lam F(\lam)=N\int_\R h(t)u^{1,\lam h}_N(t)\ud t,$$ to prove \eqre{Johlim} for derivatives and to deduce \eqre{Johlim} by dominated convergence. We use \eqre{1731418h}, namely the convergence  $$N\bet\int\ff{z-t}\lf(u_N^{1,h}(t)-\one_{|t|\le 2}\f{\sqrt{4-t^2}}{2\pi}\ri)\ud t\lto \int \ff{z-t}\ud \nu(t),$$ in Fourier transform manipulations, precisely via  the formula $$\int_0^{\infty}\hat{\del}_N(\xi)\mre^{\ii\xi z}\ud \xi=\ii\int\ff{z-t}\del_N(t)\ud t \qquad (\Im z >0)$$ with $\del_N(t):=N\bet(u_N^{1,h}(t)-\one_{|t|\le 2}\f{\sqrt{4-t^2}}{2\pi})$.

 \section{Appendix: Poisson limit for point processes}
 Let   $\X$ be a locally compact Polish space and   $\mu$ be a Radon measure on $\X$. We consider an exchangeable  random vector $(\lam_1, \ld, \lam_N)$ taking values on $\X$ implicitly depending on $N$, with density $\rho^{(N)}$ with respect to $\mu^{\otimes N}$.
  We define, for $1\le k\le N$, the $k$-th correlation function on $\X^k$ by the formula
 $$R^{(N)}_k(x_{1}, \ld, x_k):=\f{N!}{(N-k)!}\int_{(x_{k+1},\ld, x_N)\in \X^{N-k} }\rho^{(N)}(x_1, \ld, x_N)\ud \mu^{\otimes N-k}(x_{k+1},\ld, x_N).$$
 
    \beg{propo}\la{prop:convPoisson}Suppose that there is $\tta\ge 0$ independent of $N$ \st the correlation functions $R_k^{(N)}$ satisfy:\bgt\ite[a)] For each $k\ge 1$, on $\X^k$, we have the pointwise convergence \be\la{861519h}R^{(N)}_k(x_{1}, \ld, x_k)\;\Ninf \;\tta^k,\ee    \ite[b)] For each compact $\mc{K}\subset\X$, there is $\Tta_{\mc{K}}$ \st for  all $k,N$, on $\mc{K}^k$, we have   \be\la{861519h1}\one_{k\le N}R^{(N)}_k(x_{1}, \ld, x_k) \;\le\;\Tta_{\mc{K}}^k\ee
  \ent 
  Then the point process $\sum_{i=1}^N\del_{\lam_i}$ converges  in distribution to a Poisson point process with intensity $\tta\ud\mu$ as $N\to\infty$.
  \en{propo}
  
    \bpr Note that the Poisson point process $M$ with intensity $\tta\ud\mu$ is characterized, among  random random Radon measures on $\X$, by the fact that  for any compactly supported continuous function $f$ on $\X$, we have $$\E \mre^{\lan M, f\ran}=\exp \lf(\tta\int(\mre^{f(x)}-1)\ud\mu(x)\ri).$$
  So let us fix $f$ a compactly supported continuous function on $\X$.
 Then, with the convention $R_0^{(N)}=1$,  
 \beq 
\E \mre^{\sum_{i=1}^Nf(\lam_i)} 
 &=&\E\prod_{i=1}^N(1+(\mre^{f(\lam_i)}-1))\\
 &=&\sum_{P\subset\{1, \ld, n\}}\E \prod_{i\in P}(\mre^{f(\lam_i)}-1)\\
  &=&\sum_{k=0}^N\binom{N}{k}\E \prod_{i=1}^k(\mre^{f(\lam_i)}-1)\\
  &=&\sum_{k=0}^N\ff{k!}\int \prod_{i=1}^k(\mre^{f(x_i)}-1) R^{(N)}_k(x_{1}, \ld, x_k)\ud\mu^{\otimes k}(x_1, \ld, x_k)
 \eeq 
 This proves the proposition.\epr

    \begin{thebibliography}{10}
    
    \bibitem{AllezGuionnetBouchaud} R. Allez, J.-P. Bouchaud, A. Guionnet \emph{Invariant Beta ensembles and the Gauss-Wigner crossover}
 Physical review letters 109 (9), 094--102.
 
 \bibitem{Allez2} R. Allez, J.-P. Bouchaud, S.N. Majumdar, P. Vivo \emph{Invariant $\beta$-Wishart ensembles, crossover densities and asymptotic corrections to the Marcenko-Pastur law}
Journal of Physics A: Mathematical and Theoretical 46 (1), 015001

\bibitem{AllezGuionnet} R. Allez,  A Guionnet \emph{A diffusive model for invariant $\beta$-ensembles}
  Electron. J. Probab 18 (62), 1--30.

\bibitem{AD1} R. Allez, L. Dumaz \emph{From sine kernel to Poisson statistics} Elec. J. Probab. (2014), Vol. 19,  1--25.

\bibitem{AD2} R. Allez, L. Dumaz \emph{Tracy--Widom at High Temperature} J. Stat. Phys. (2014), Vol. 156, No 6,   1146--1183.

 \bibitem{alice-greg-ofer}
G.~Anderson, A.~Guionnet, O.~Zeitouni \emph{An Introduction to Random Matrices}. Cambridge studies in advanced mathematics, {118} (2009).

\bibitem{bao} Z. Bao, Z. Su \emph{Local Semicircle law and Gaussian fluctuation for Hermite $\beta$-ensemble}, arXiv.
   
   \bibitem{BAG3}
G.~Ben Arous, A.~Guionnet.
\emph{Large deviations for  Wigner's law and Voiculescu's non commutative entropy},
{Probab. Theory Related Fields},
{\bf 108}, (1997) 517--542.

  \bibitem{BLM} S. Boucheron, G. Lugosi, P. Massart \emph{Concentration inequalities}, Oxford, 2013.
 
 \bibitem{BEY4} P. Bourgade  \emph{Bulk universality for one-dimensional log-gases}, proceedings of the XVIIth International Congress On Mathematical Physics, World Scientific Publishing, 2013.

\bibitem{BEY1} P. Bourgade, L. Erd\"os, H.T. Yau \emph{Bulk universality of general $\beta$-ensembles with non-convex potential}, J. of Math. Phys., special issue in honor of E. Lieb's 80th birthday, Vol. 53, 2012.

\bibitem{BEY2} P. Bourgade, L. Erd\"os, H.T. Yau \emph{Universality of general $\beta$-ensembles}, Duke Math. J., Vol. 163, no. 6, 1127--1190, 2014.

\bibitem{BEY3} P. Bourgade, L. Erd\"os, H.T. Yau \emph{Edge Universality of $\bet$-ensembles}, to appear in Comm. Math. Phys., 2013.

  %\bibitem{DZ} A.~ Dembo, O.~ Zeitouni, \emph{Large Deviation Techniques and Applications}. {Applications of Mathematics (New York)} {\bf 38}, Springer (1998).
  
  \bibitem{Durrett2010} R. Durrett \emph{Probability : Theory and examples}, Fourth edition, Cambridge Univ. Press, 2010.

  \bibitem{DumEdel} I. Dumitriu,  A. Edelman \emph{Matrix models for beta ensembles}. J. Math. Phys., 43:5830--5847, (2002).
  
  \bibitem{EnriquezMenard} N. Enriquez, L. M\'enard \emph{Asymptotic expansion of the expected spectral measure of Wigner matrices}, arXiv.

  \bibitem{Johansson} K. Johansson \emph{On fluctuations of eigenvalues of random Hermitian matrices} Duke Math. Journ. 91, no.1, 151-204 (1998).
  
  \bibitem{KillipStoicu} R. Killip, M. Stoiciu \emph{Eigenvalue statistics for CMV matrices: From Poisson to clock via random matrix ensembles},		Duke Math. J. 146 (2009), 361--399.
  
  \bibitem{LebleSerfaty} T. Lebl\'e, S. Serfaty \emph{Large Deviation Principle for Empirical Fields of Log and Riesz Gases}
arXiv:1502.02970 (2015)
  
  \bibitem{Ns06} A. Nica, R. Speicher  \emph{Lectures on the combinatorics of free probability}. London Mathematical Society Lecture Note Series, 335. Cambridge University Press, Cambridge, 2006.
  
  \bibitem{RamirezRiderVirag} J.A. Ram\'{\i}rez, B. Rider, B. Vir\'ag \emph{Beta ensembles, stochastic Airy spectrum, and a diffusion}.  J. Amer. Math. Soc. 24 (2011), 919--944 

 \bibitem{sosoe} P. Sosoe, P. Wong,  \emph{Local semicircle law in the bulk for Gaussian $\beta$-ensemble}. J. Stat. Phys. 148 (2012), no. 2, 204--232.
 
 \bibitem{sosoe2} P. Sosoe, P. Wong,  \emph{Convergence of the eigenvalue density for $\beta$-Laguerre ensembles on short scales}. Electron. J. Probab. 19 (2014), no. 34, 18 pp.
  
     \bibitem{SimonTrace} {B. Simon} \emph{Trace ideals and their applications}, {Mathematical Surveys and Monographs}, {120},  {Second Edition}, AMS, 2005.
     
     \bibitem{TrinhTomoyuki} K.D. Trinh, S.  Tomoyuki \emph{The mean spectral measures of random Jacobi matrices related to Gaussian beta ensembles}. arXiv

\bibitem{ValkoVirag} B. Valk\'o, B. Vir\'ag \emph{Continuum limits of random matrices and the Brownian carousel.}   Invent math (2009) 177: 463-508.

 \en{thebibliography}
\en{document}